# SHARP ADAPTIVE ESTIMATION OF THE DRIFT FUNCTION FOR ERGODIC DIFFUSIONS

By Arnak Dalalyan

*Université Paris VI*

The global estimation problem of the drift function is considered for a large class of ergodic diffusion processes. The unknown drift $S(\cdot)$ is supposed to belong to a nonparametric class of smooth functions of order $k \geq 1$, but the value of $k$ is not known to the statistician. A fully data-driven procedure of estimating the drift function is proposed, using the estimated risk minimization method. The sharp adaptivity of this procedure is proven up to an optimal constant, when the quality of the estimation is measured by the integrated squared error weighted by the square of the invariant density.

## 1. Introduction.

1.1. *The problem.* In this paper we consider the statistical problem of estimating the drift function of a diffusion process $X$, given as the solution of the stochastic differential equation

$$(1) \qquad dX_t = S(X_t)\,dt + \sigma(X_t)\,dW_t, \qquad X_0 = \xi, t \geq 0,$$

where $W$ is a standard Brownian motion and the initial value $\xi$ is a random variable independent of $W$. We assume that a continuous record of observations $X^T = (X_t, 0 \leq t \leq T)$ is available. The goal is to estimate the function $S(\cdot)$, which is commonly referred to as the drift function and is interpreted as the instantaneous mean of the process $X$.

In our setup, the diffusion coefficient $\sigma^2$ is identifiable using the quadratic variation of the semi-martingale $X$. Therefore, the problem of its estimation is not interesting from the viewpoint of asymptotic statistics. In the sequel, we suppose that $\sigma(\cdot)^2$ is a known function satisfying some boundedness and smoothness conditions.

During the past decade statistical inference for continuous time Markov processes has been widely developed due to its numerous applications, namely, in mathematical finance and econometrics. In fact, the diffusion processes and their extensions, such as jump-diffusions and the solutions of stochastic differential equations driven by Lévy processes, are often used to model the evolution of asset prices and derivative securities.

Estimation problems based on both continuous time and discretely sampled observations have been considered in the statistical literature. The first one is more interesting from the point of view of financial econometrics (see [1, 25]), since in finance, even if the underlying process is time continuous, only its values at a finite number of points are available.

However, the theoretical development of statistical inference with a continuous time record of observations turns out to be technically simpler than inference based on discretely observed data. It permits one, therefore, to go further in the statistical analysis of the model and to answer some questions open up to now for discretely observed diffusions.

Note also that the continuous-time model can be considered as the limit of a discrete-time model when the step of discretization goes to zero (see [30]). Therefore, if the available data is "dense enough" with respect to the observation time, the asymptotic behavior of estimation procedures, in practice, may be close to the asymptotic behavior proven theoretically for continuous-time observations. Thus, the knowledge of the best estimator based on continuous-time data is of practical interest as well.

The purpose of the present paper is to estimate the drift function globally, that is, at any point $x \in \mathbb{R}$. We consider the case of ergodic diffusions, which means that the Markov process $X$ admits an invariant measure. Let $f_S$ denote the density of this invariant measure with respect to the Lebesgue measure on $\mathbb{R}$ (cf. [18], Chapter 4, Section 18, for more details). To quantify the performance of an estimator $S_T(\cdot) = S_T(\cdot, X^T)$ of the drift $S(\cdot)$, we use the weighted $L^2$-risk

$$(2) \qquad R_T(S_T, S) = \int_{\mathbb{R}} \mathbf{E}_S[(S_T(x) - S(x))^2] f_S^2(x) \, dx,$$

where $\mathbf{E}_S$ is the expectation with respect to the law $\mathbf{P}_S$ of $X$ defined by (1). We call an estimating procedure *adaptive* if its realization does not require any a priori information on the estimated function. The only information that we may (and should) use is that contained in the observations. We call an estimating procedure *minimax sharp adaptive*, or simply *sharp adaptive*, if its minimax risk converges with the best possible rate to the best possible constant.

The main focus of this paper is constructing an estimating procedure which is minimax sharp adaptive with respect to the risk (2), when $T \to \infty$.

The estimator of the drift we propose enjoys the following properties:

– It is fully data-driven; particularly, it does not depend on the smoothness of the estimated function.
– In the case when the observed path is generated by a stationary diffusion, the estimator is minimax rate-optimal over the Sobolev balls on any compact interval. The quality of estimation is then measured by the mean integrated squared error (MISE).
– Still in the case of an ergodic diffusion, if the risk function is defined by (2), the estimator is asymptotically sharp adaptive over a large scale of Sobolev balls: it attains not only the optimal rate of convergence, but also the best possible constant. Moreover, the accuracy of the adaptive procedure is asymptotically as good as the accuracy of the best possible nonadaptive procedure.

1.2. *Adaptive estimation.* The first result concerning minimax sharp adaptivity in nonparametric curve estimation is due to Efromovich and Pinsker [11]. It was extended by Golubev [19] and Golubev and Nussbaum [22] for nonparametric regression and by Efromovich [9] and Golubev [21] for density estimation from i.i.d. data. Similar results have been obtained in some other contexts as well (we refer to Chapter 7.4 of [10] for a comprehensive discussion), but they all deal with either independent or Gaussian observations. Thus, the main difference of our study is that the observations we have at our disposal are neither independent nor Gaussian. Moreover, as follows from heuristics presented in Section 4.1, our model exhibits heteroscedastic structure.

1.3. *Estimation for diffusions.* For a complete review of parametric and nonparametric methods for diffusion processes, we refer to [14], [24] and [27]. There are a number of papers devoted to the estimation of the drift in the case when the parameters (such as smoothness, Lipschitz constant) describing the nonparametric class are known and when continuous-time observations are available. Banon [3] proved the consistency in probability of kernel-type estimators, and Pham [31] obtained the rate of convergence in the same setup. These results have been extended by van Zanten [35], Galtchouk and Pergamenshchikov [17] and Kutoyants [27]. Pinsker's constant in this problem is obtained in [6]. More recently, an approach making use of a random rate of convergence is developed in [8].

The adaptive estimation of the drift at a fixed point based on continuous-time observations has been studied by Spokoiny [34], who has applied the Lepskii method (see [28]) to the locally linear smoothers in order to construct an adaptive rate-optimal procedure. For discretely sampled diffusions, a rate-optimal adaptive procedure for estimating the drift function has been proposed by Hoffmann [23].

In the problem of estimation and validation of the model with discretely sampled high frequency observations, which is frequently used in mathematical finance, recent progress has been achieved by Fan and Zhang [15] and Aït-Sahalia and Mykland [2].

Note also that the diffusion processes can be regarded as continuous-time analogues of autoregressive processes. A theoretical result establishing the connection between these two models has been proven by Milstein and Nussbaum [29].

The paper is divided into five sections. The description of the adaptive procedure is given in Section 2. In Section 3 the assumptions on the model, as well as the main result describing the asymptotic behavior of the procedure, are formulated. Some comments related to this result and its proof are given in Sections 4 and 5, respectively.

## 2. Construction of the adaptive estimator.

The construction of the adaptive estimator proposed in this paper relies hardly on the papers [4] and [21]. It can be divided into three steps according to the following scheme. First, we present an estimator of the drift function involving two kernel-type functions and two bandwidths. This estimator is rate-optimal if the orders of the bandwidths are chosen in a correct way. Second, we derive an asymptotically exact upper bound of the risk of this estimator. The minimization of the maximum of this risk bound over the Sobolev ball of smoothness $k$ and radius $R$ provides the explicit forms of the optimal kernel and the optimal bandwidth, depending on $k$ and $R$. The last step is to substitute some "good" data-driven approximation for these parameters.

2.1. *The nonadaptive estimator.* Let $K(\cdot), Q(\cdot) \in L^2(\mathbb{R})$ be two positive $k$-times ($k \geq 1$) continuously differentiable symmetric functions such that $\int K = \int Q = 1$, and let $\alpha = \alpha_T$ and $\nu = \nu_T$ be two positive functions of $T$ decreasing to zero. According to the Girsanov formula, for two drift functions $S$ and $S_0$, the log-likelihood $\log(\frac{d\mathbf{P}_S}{d\mathbf{P}_{S_0}}(X^T))$ in the model (1) is given by

$$\Lambda_T(S, S_0, X^T) = \int_0^T \frac{S(X_t) - S_0(X_t)}{\sigma^2(X_t)}\, dX_t - \frac{1}{2} \int_0^T \frac{S^2(X_t) - S_0^2(X_t)}{\sigma^2(X_t)}\, dt.$$

We have supposed in the above formula that the law of the initial value does not depend on $S$. A widely used idea for constructing nonparametric estimators is to find a local (around the state $x$) approximation $\tilde{\Lambda}(\theta, X^T, x)$ of $\Lambda_T(S, S_0, X^T)$ depending only on a finite-dimensional parameter $\theta$ and to define the estimator as the value of the parameter $\theta$ maximizing $\tilde{\Lambda}(\theta, X^T, x)$. We use the following "local constant" approximation of the log-likelihood:

$$\tilde{\Lambda}(\theta, X^T, x) = \frac{\theta}{\alpha_T \sigma^2(x)} \int_0^T K\left(\frac{x - X_t}{\alpha_T}\right) dX_t - \frac{\theta^2}{2\nu_T \sigma^2(x)} \int_0^T Q\left(\frac{x - X_t}{\nu_T}\right) dt$$

(in this expression, the terms not depending on $\theta$ are dropped, since they have no influence on the definition of the MLE). It is evident that the maximum of this expression is attained by

$$(3) \qquad \theta_T(x) = \frac{1}{\alpha_T} \int_0^T K\left(\frac{x - X_t}{\alpha_T}\right) dX_t \left[\frac{1}{\nu_T} \int_0^T Q\left(\frac{x - X_t}{\nu_T}\right) dt\right]^{-1},$$

provided that the denominator is different from 0. A similar algorithm but with local linear smoothers is used in [34]. As it is explained in Section 4.4 of [6], for symmetric functions $K$ and $Q$, the asymptotic properties of the estimators defined via the local constant and the local linear smoothers coincide. That is why we restrict ourselves to the local constant approximation.

One drawback of the estimator (3) is that it is not defined when the denominator is null. Different approaches for overcoming this problem have been proposed. Banon [3] has suggested increasing artificially the denominator by a deterministic term $\epsilon/(T\nu_T)$, which asymptotically vanishes but allows the denominator to stay positive. We adopt in this paper this approach, but with a more careful choice of the term to be added. Note that, in the context of nonparametric regression, a similar approach is developed in [13].

The second drawback of the estimator (3) is the presence of the stochastic integral. To explain why this integral is undesirable, let us recall that our final goal is to choose all the parameters and, in particular, the bandwidth $\alpha_T$ in a data-dependent way. If we replace $\alpha$ by an approximation depending on the observed path $(X_t, 0 \le t \le T)$, we obtain an anticipative stochastic integral. The manipulation of such integrals is technically more difficult than the manipulation of the Riemann integrals.

In order to replace the stochastic integral by a Riemann integral, we apply the Itô formula to the primitive of the function $\alpha^{-1} K((x - \cdot)/\alpha)$ and to the semi-martingale $X$:

$$\int_{X_0}^{X_T} K\left(\frac{x - y}{\alpha}\right) dy = \int_0^T K\left(\frac{x - X_t}{\alpha}\right) dX_t - \frac{1}{2\alpha} \int_0^T K'\left(\frac{x - X_t}{\alpha}\right) \sigma^2(X_t) dt.$$

We show that, in the ergodic case, the term $\int_{X_0}^{X_T} K((x - y)/\alpha) \, dy$ is asymptotically negligible with respect to the other terms. Therefore, the stochastic integral $\int_0^T K((x - X_t)/\alpha) \, dX_t$ can be approximated by $(2\alpha)^{-1} \int_0^T K'((x - X_t)/\alpha)\sigma^2(X_t) \, dt$.

According to these considerations, we modify the estimator (1) as

$$(4) \qquad \hat{\theta}_T(x) = \frac{(1/\alpha^2) \int_0^T K'((x - X_t)/\alpha)\sigma^2(X_t) \, dt}{(2/\nu) \int_0^T Q((x - X_t)/\nu) \, dt + (2\varepsilon/\nu)e^{-\ell_T |x|}},$$

where $\varepsilon = \varepsilon_T = e^{\sqrt{\log T}}$ and $\ell_T = (\log T)^{-1}$. It is proved in [27] that if the unknown drift function is $k$-times continuously differentiable, then the bandwidths $\alpha_T = T^{-1/(2k+1)}$ and $\nu_T = T^{-1/2}$ lead to a locally and globally rate-optimal estimator $\hat{\theta}_T(x)$. The rate of convergence of this estimator is then $T^{-k/(2k+1)}$. This choice of the bandwidth $\alpha_T$ is clearly nonadaptive, since it depends on the unknown parameter $k$.

In the case of an ergodic diffusion, one can arrive at the same estimator using the well-known formula

$$(5) \qquad (\sigma^2(x)f_S(x))' = 2S(x)f_S(x),$$

where $f_S$ is the invariant density of the process $X$. Using the occupations time formula and the martingale representation of the local time, one can check that $(T\alpha^2)^{-1}\int_0^T K'((x-X_t)/\alpha)\sigma^2(X_t)\,dt$ is a consistent estimator of $(\sigma^2(x)f_S(x))'$. Likewise, $2(T\nu)^{-1}\int_0^T Q((x-X_t)/\nu)\,dt$ is a consistent estimator of $2f_S(x)$. It is now quite natural to define the estimator of $S(x)$ as the quotient of these two estimators.

2.2. *The minimax sharp adaptive estimator.* To simplify the exposition, we suppose in this section that the diffusion coefficient $\sigma(\cdot)$ is identically equal to one. For any function $h \in L^2(\mathbb{R})$, let us denote by $\varphi_h(\cdot)$ the Fourier transform of $h(\cdot)$ defined as $\varphi_h(\lambda) = \int_{\mathbb{R}} e^{i\lambda x}h(x)\,dx$. To avoid double subscripts, we write $\varphi_f$ instead of $\varphi_{f_S}$. Recall that, for any estimator $S_T(\cdot) = S_T(\cdot, X^T)$ of the drift function $S(\cdot)$, we have defined

$$R_T(S_T, S) = \int_{\mathbb{R}} \mathbf{E}_S[(S_T(x) - S(x))^2]f_S^2(x)\,dx.$$

Some heuristic explanations of this choice of the risk function are presented in Section 4.1. It is proven in [6] that, in order that the estimator (4) be asymptotically minimax over a properly chosen Sobolev ball $\Sigma(k, R)$ ($k$ is the order of smoothness and $R$ is the radius), one should choose the kernels and the bandwidths as

$$
(6) \qquad
\begin{aligned}
\alpha_T^* &= \left(\frac{4k}{\pi RT(k+1)(2k+1)}\right)^{1/(2k+1)}, \\
K^*(x) &= \frac{1}{\pi}\int_0^1 (1 - u^{k+\rho_T})\cos(ux)\,du;
\end{aligned}
$$

$\nu_T = T^{-1/2}$ and $Q(x)$ is any positive, differentiable, symmetric function with support in $[-1, 1]$ and $\int Q(x)\,dx = 1$. In equality (6), we used the notation $\rho_T = 1/\log\log(1 + T)$. The estimator (4) defined by such a bandwidth and kernel will be denoted by $S_T^*(\cdot)$. Note here that the Fourier transform of the kernel $K^*$ is $\varphi_{K^*}(\lambda) = (1 - |\lambda|^{k+\rho_T})_+$. The exact asymptotic behavior of the maximum over $\Sigma(k, R)$ of the risk of this estimator is $T^{-2k/(2k+1)}P(k, R)$ (see [6], Theorem 4 and Definition 2), where $P(k, R)$ is Pinsker's constant [32]. Moreover, the following asymptotic relation holds:

$$R_T(S_T^*, S) \leq \frac{\Delta_T(\alpha, \varphi_{K^*}, \varphi_f)(1 + o_T(1))}{2\pi T},$$

where $o_T(1)$ is a term tending to zero uniformly in $S$ and the functional $\Delta_T$ is defined by

$$\Delta_T(\alpha, h, \varphi_f) = T \int_{\mathbb{R}} |\lambda(1 - h(\alpha\lambda))\varphi_f(\lambda)|^2 \, d\lambda + 4 \int_{\mathbb{R}} |h(\alpha\lambda)|^2 \, d\lambda.$$

Since for known $k$ the optimal kernel is given by (6), it will be natural to select the adaptive kernel among the functions $\{K_\beta(x) = \pi^{-1} \int_0^1 (1 - u^\beta) \times \cos(ux) \, du | \beta > 0\}$ in a data-driven way. Thus, it suffices to give a good adaptive choice of the real parameters $\alpha$ and $\beta$ in order to obtain an adaptive estimator of $S$. The values of these parameters that are of interest for us are those minimizing the risk $R_T(S_T, S)$ or, equivalently, $\Delta_T(\alpha, h_\beta, \varphi_f)$, where

$$h_\beta(\lambda) = (1 - |\lambda|^\beta)_+.$$

The minimizers of $\Delta_T$ depend obviously on the unknown function $S$, so they cannot be used in an estimation procedure. A standard method for overcoming this difficulty is to estimate $\Delta_T(\alpha, h, \varphi_f)$ by a data dependent functional $l_T(\alpha, h)$ that does not involve the function $S$. Then the minimizers of the latter functional might be chosen as parameters for the adaptive procedure. Perhaps the most straightforward idea for estimating $\Delta_T(\alpha, h, \varphi_f)$ is to utilize the plug-in estimator $\Delta_T(\alpha, h, \hat{\varphi}_T)$, $\hat{\varphi}_T$ being the empirical characteristic function

$$\hat{\varphi}_T(\lambda) = \frac{1}{T} \int_0^T e^{i\lambda X_t} \, dt.$$

But it is well known that the plug-in estimators of quadratic functionals have a large bias (cf., e.g., [12]). That is why a smarter solution consists in applying the plug-in method to $\Delta_T$ considered as a linear functional of $|\varphi_f(\cdot)|^2$. According to Lemma 1, a good estimate of $|\varphi_f(\lambda)|^2$ is $|\hat{\varphi}_T(\lambda)|^2 - 4/(T\lambda^2)$. On the other hand, the minimization of $\Delta_T(\alpha, h_\beta, \varphi_f)$ w.r.t. parameters $\alpha$ and $\beta$ is obviously equivalent to the minimization of

$$\tilde{\Delta}_T(\alpha, h_\beta, \varphi_f) = T \int_{\mathbb{R}} \lambda^2 (h_\beta^2(\alpha\lambda) - 2h_\beta(\alpha\lambda)) |\varphi_f(\lambda)|^2 \, d\lambda + 4 \int_{\mathbb{R}} |h_\beta(\alpha\lambda)|^2 \, d\lambda,$$

since it is just $\Delta_T(\alpha, h_\beta, \varphi_f) - T \int_{\mathbb{R}} \lambda^2 |\varphi_f(\lambda)|^2 \, d\lambda$. For this reason, we define the functional

$$
\begin{aligned}
l_T(h) = {}& T \int_{\mathbb{R}} \lambda^2 (h^2(\lambda) - 2h(\lambda)) |\hat{\varphi}_T(\lambda)|^2 \, d\lambda - 4 \int_{\mathbb{R}} (h^2(\lambda) - 2h(\lambda)) \, d\lambda \\
& + 4 \int_{\mathbb{R}} |h(\lambda)|^2 \, d\lambda \\
= {}& T \int_{\mathbb{R}} \lambda^2 (h^2(\lambda) - 2h(\lambda)) |\hat{\varphi}_T(\lambda)|^2 \, d\lambda + 8 \int_{\mathbb{R}} h(\lambda) \, d\lambda.
\end{aligned}
\tag{7}
$$

This functional depends on the observed path via the empirical characteristic function $\hat{\varphi}_T$. To obtain the adaptive kernel $K_\beta$ and the adaptive bandwidth $\alpha$, one should minimize the expression $l_T(h)$ over a suitably chosen subset $\mathcal{H}_T^N$ of the set

$$\mathcal{H}_T = \{h : x \mapsto h_\beta(\alpha x) = (1 - |\alpha x|^\beta)_+ | \alpha \in [T^{-1/3}, (\log T)^{-1}], \beta \geq 1\},$$

such that $\#\mathcal{H}_T^N = N$. The subset $\mathcal{H}_T^N$ is defined as in [4]. For any pair of positive integers $i$ and $j$, let us denote

$$(8) \qquad \alpha_i = \left(1 + \frac{1}{\log T}\right)^{-i} \quad \text{and} \quad \beta_j = \left(1 - \frac{j}{\log T}\right)^{-1}.$$

The finite subset $\mathcal{H}_T^N$ of $\mathcal{H}_T$ is defined as

$$\mathcal{H}_T^N = \{h : x \mapsto (1 - |\alpha_i x|^{\beta_j})_+ | \alpha_i \in [T^{-1/3}, (\log T)^{-1}], j = 1, \ldots, \lfloor \log T \rfloor\},$$

where $\lfloor a \rfloor$ denotes the largest integer strictly smaller than the real number $a$. It is evident that the cardinality of $\mathcal{H}_T^N$ is less than $(\log T)^3$. From now on, we denote the $N$ elements of this set by $h_1, h_2, \ldots, h_N$. Thus, to construct the adaptive estimator, the functional $l_T$ is maximized over a set of cardinality not exceeding $(\log T)^3$. Let us now summarize the method.

2.3. *Brief description of the procedure.* We start by computing the values $\alpha_i$ and $\beta_j$ according to (8). Then we determine the function $\tilde{h}_T \in \mathcal{H}_T^N$ such that $l_T(\tilde{h}_T) = \min_{h \in \mathcal{H}_T^N} l_T(h)$. If the function satisfying the latter equality is not unique, we denote by $\tilde{h}_T$ one of them. Next we apply the inverse Fourier transform to $\tilde{h}_T$ in order to define the kernel

$$\tilde{K}_T(x) = \frac{1}{2\pi} \int_\mathbb{R} \tilde{h}_T(\lambda) \cos(\lambda x) \, d\lambda.$$

This form of the kernel comprises the bandwidth since $\tilde{h}_T(\lambda) = h_{\tilde{\beta}_T}(\tilde{\alpha}_T \lambda)$, where $\tilde{\alpha}_T$ and $\tilde{\beta}_T$ are the values of $\alpha_i$ and $\beta_j$ corresponding to $\tilde{h}_T$. Further, we choose another kernel function $Q(\cdot)$ which is positive, differentiable, symmetric, supported in $[-1, 1]$ and with integral equal to one. Finally, we set $\varepsilon_T = e^{\sqrt{\log T}}$, $\ell_T = 1/\log T$ and define the estimator

$$\tilde{S}_T(x) = \frac{\int_0^T \tilde{K}_T'(x - X_t) \, dt}{2\sqrt{T} \int_0^T Q(\sqrt{T}(x - X_t)) \, dt + 2\sqrt{T} \varepsilon_T e^{-\ell_T |x|}}.$$

Note that the function $\tilde{K}_T(\cdot)$ is differentiable, since $\min \beta_j > 1$.

**3. Assumptions and main results.** We introduce four conditions playing an important role throughout this paper. They ensure the existence of the observed diffusion process as a solution of (1) and provide us with some technical tools permitting us to deal with this process. Before stating these conditions, we need some additional notation.

Recall that the solution of the stochastic differential equation (1) is a strong Markov process. We denote by $P_t(S, x, A)$ the transition probability corresponding to the instant $t$, that is,

$$P_t(S, x, A) = \mathbf{P}_S(X_t \in A | X_0 = x) \qquad \forall\, x \in \mathbb{R}, \forall\, A \in \mathscr{B}(\mathbb{R}).$$

Here $\mathbf{P}_S$ denotes the probability measure on $(C(\mathbb{R}), \mathscr{B}_{C(\mathbb{R})})$ induced by the process (1). For every $x \in \mathbb{R}$ and $t \geq 0$, the probability measure $P_t(S, x, \cdot)$ is absolutely continuous with respect to the Lebesgue measure. The corresponding density will be denoted by $p_t(S, x, y)$, so that, for any integrable function $g(\cdot)$, we have

$$(9) \qquad \mathbf{E}_S[g(X_t)|\mathscr{F}_s] = \int_{\mathbb{R}} g(y)\, p_{t-s}(S, X_s, y)\, dy.$$

Let $k$ be a strictly positive integer. Denote by $\Sigma(k)$ the set of all functions satisfying the following conditions:

C1. The function $S$ is $k$-times continuously differentiable in the whole real line and $\limsup_{|x| \to \infty} S(x)\, \mathrm{sgn}\, x < 0$.

C2. There exist positive numbers $C$ and $\nu$ such that $|S^{(k)}(x)| \leq C(1 + |x|^\nu)$, $\forall\, x \in \mathbb{R}$.

The problem we consider is the following: we know that $x^T$ is a sample path of the process $X^T$ given by (1) with a drift function $S \in \Sigma = \bigcup_{k \geq 1} \Sigma(k)$ and we want to estimate the function $S(\cdot)$. To obtain minimax results, we consider the local setting. For any function $S_0 \in \Sigma(k)$ and for all $\delta > 0$, we define the neighborhoods $V_\delta(S_0) = \{S \in \Sigma | \sup_{x \in \mathbb{R}} |S(x) - S_0(x)| \leq \delta\}$ and

$$\widetilde{V}_\delta(S_0) = \left\{ S \in \Sigma(k) \Big| \sup_{x \in \mathbb{R}} |S^{(i)}(x) - S_0^{(i)}(x)| \leq \delta, i = 0, 1, \ldots, k-1 \right\}.$$

The center of localization $S_0(\cdot)$ is assumed to fulfill the following additional assumptions:

C3. There exist a positive number $\kappa$ and a $q > 1$ such that the quantity $\sup_{t > \kappa} \mathbf{E}_{S_0}[\sup_{y \in \mathbb{R}} p_\kappa^q(S_0, X_{t-\kappa}, y)]$ is finite.

C4. Let $\varphi_0(\cdot)$ be the Fourier transform of the invariant density $f_{S_0}(\cdot)$. There exists $\tau > 0$ such that $\int_{\mathbb{R}} |\lambda|^{2k+2+\tau} |\varphi_0(\lambda)|^2\, d\lambda < \infty$.

Conditions C1–C4 need perhaps some comments. The first one ensures the ergodicity (see [18]) of the solution of the stochastic differential equation (1).

This condition entails also the exponential smallness of the tails of $f_S(\cdot)$. The second condition guarantees the square integrability of the functions $f_S^{(i)}(\cdot)$, for every $i = 0, 1, \ldots, k$.

Condition C3 is a technical one and can be considered a mixing property of the underlying diffusion process. It can be viewed as a weakened version of the condition $G_2(s, \alpha)$ from [3]. Some sufficient conditions for C3 are given in Section 4.4.

Finally, C4 means that the central function $S_0(\cdot)$ is a little bit smoother than the other functions of the neighborhood. For example, if $S_0(\cdot)$ is $(k+1)$-times differentiable and $S_0^{(k+1)}(\cdot)$ increases at most polynomially, then C4 is satisfied with $\tau = 2$.

We define now the Sobolev balls; in our setup they also are weighted by the square of the invariant density. Let us denote

$$\widetilde{\Sigma}_\delta(k, R, S_0) = \left\{ S \in \widetilde{V}_\delta(S_0) \Big| \int_{\mathbb{R}} [(S - S_0)^{(k)}(x)]^2 f_S^2(x) \, dx \le R \right\}.$$

To simplify the notation, we write $\widetilde{\Sigma}_\delta$ instead of $\widetilde{\Sigma}_\delta(k, R, S_0)$.

We state now the main theorem of this work describing the asymptotic behavior of the estimator $\tilde{S}_T$ constructed in the previous section.

THEOREM 1.   *Let $S_0$ satisfy assumptions* C1–C4 *and let the risk $R_T(\cdot, \cdot)$ be defined by* (2). *If the initial condition $\xi$ follows the invariant law, then*

$$\limsup_{\delta \to 0} \limsup_{T \to \infty} \sup_{S \in \widetilde{\Sigma}_\delta} T^{2k/(2k+1)} R_T(\tilde{S}_T, S) = P(k, R),$$

*where $P(k, R) = (2k + 1)(\frac{k}{\pi(k+1)(2k+1)})^{2k/(2k+1)} R^{1/(2k+1)}$ is Pinsker's constant (cf. [32]).*

Note that this asymptotic bound cannot be improved, since it coincides with the lower bound obtained in [6], Theorem 3. Hence, the adaptive estimator $\tilde{S}_T$ behaves asymptotically as well as the best possible nonadaptive estimator, provided that the specific form (2) of the risk function is used.

*Consequence.*   As an immediate consequence of the above theorem, one obtains the rate-optimality of the estimator $\tilde{S}_T$ when the error of estimation is quantified using the MISE over a compact set $\mathcal{K} \subset \mathbb{R}$. That is, for sufficiently small values of $\delta$, we have

$$\limsup_{T \to \infty} \sup_{S \in \widetilde{\Sigma}_\delta} T^{-2k/(2k+1)} \mathbf{E}_S \int_{\mathcal{K}} (\tilde{S}_T(x) - S(x))^2 \, dx \le C < \infty.$$

The reason for this is the uniform in $S \in \widetilde{\Sigma}_\delta$ boundedness of the functions $f_S$ and $f_S^{-1}$ on any compact set $\mathcal{K}$.

## 4. Remarks and extensions.

4.1. *The weight function.* The choice of the weight function in the risk definition is mainly motivated by the weak equivalence of experiments. In fact, as is explained in detail in [6], the Gaussian white noise experiment having almost the same statistical properties as our model is

$$(10) \qquad dY_t = S(t)\,dt + [T f_0(t)]^{-1/2}\,dB_t, \qquad t \in \mathbb{R},$$

where $B_t$ is a two-sided standard Brownian motion and $f_0 = f_{S_0}$. On the other hand, according to Golubev [20], the asymptotically optimal lower bound of the maximum of MISE (over the Sobolev ball of smoothness $k$ and radius $R$) in the model

$$dZ_t = \theta(t)\,dt + \varepsilon I^{-1/2}(\theta_0, t)\,dB_t, \qquad t \in \mathcal{I},$$

is equal to $\varepsilon^{4k/(2k+1)} P(k, R) [\int_{\mathcal{I}} I^{-1}(\theta_0, t)\,dt]^{2k/(2k+1)}$. Our aim is to find a normalization of the MISE via a weight function such that the resulting limit of the minimax risk does not depend on the central function $\theta_0$. This would hold if the integral of the Fisher information $I^{-1}(\theta_0, \cdot)$ were independent of $\theta_0$. Obviously, this is not the case for (10). In order to obtain a model enjoying the desired property, we transform (10) by multiplying it by $f_0$. We get

$$(11) \qquad d\tilde{Y}_t = S(t) f_0(t)\,dt + \sqrt{T^{-1} f_0(t)}\,dB_t, \qquad t \in \mathbb{R}.$$

The integral of the inverse of the Fisher information associated with the last model is one, since $f_0$ is a probability density. On the other hand, since the function $f_S$ is a regular functional of $S$, it can be estimated with more precision than the function $S$. At a heuristic level, this is the reason why estimating the function $S f_S$ in $L^2$ is equivalent to estimating $S$ in $L^2$ with the weight function $f_S^2$.

Note also that the use of a weight function for estimating $S$ over the whole real line is unavoidable, otherwise the risk of estimation will explode. Moreover, any deterministic weight function has to depend on the unknown function $S$ (or, at least, on an upper estimate of $S$). Indeed, if we observe a path $X^T$, it contains no information about the values of $S$ that are outside of the interval $[x_*, x^*]$, where $x_* = \min_{t \in [0,T]} X_t$ and $x^* = \max_{t \in [0,T]} X_t$. Thus, the error of estimating $S$ at a point $x \notin [x_*, x^*]$ is large when $S(x)$ is large. Consequently, in order that the integral $\int_{\mathbb{R}} (S_T - S)^2 q_S$ be finite, the weight function $q_S$ should be small when $S$ is large.

4.2. *The case of a general diffusion coefficient.* Let us consider the case where $\sigma$ is an arbitrary positive function such that $\sigma^2 + \sigma^{-2}$ is bounded by some polynomial function. It is also assumed that $\sigma$ is $(k+1)$-times differentiable and the condition C1 is replaced by $\limsup_{|x| \to \infty} S(x) \operatorname{sgn} x / \sigma^2(x) < 0$.

In this case, the functional $\Delta_T = \Delta_T(\alpha, h, \varphi_{\sigma^2 f}, \|\sigma\|^2_{L^2(f)})$ has the form

$$T \int_{\mathbb{R}} |\lambda(1 - h(\alpha\lambda)) \varphi_{\sigma^2 f}(\lambda)|^2 \, d\lambda + 4 \int_{\mathbb{R}} \sigma^2(x) f_S(x) \, dx \int_{\mathbb{R}} |h(\alpha\lambda)|^2 \, d\lambda.$$

It follows from this expression that, in order to construct an estimator $l_T$ of $\Delta_T$, one has to estimate not only the square of the Fourier transform $|\varphi_{\sigma^2 f}(\lambda)|$, but also the term $\|\sigma\|^2_{L^2(f)} = \int_{\mathbb{R}} \sigma^2(x) f_S(x) \, dx$. Fortunately, this latter quantity is just a linear functional of $f_S$ and therefore can be estimated with a parametric rate $T^{-1/2}$. Let us denote $\hat{\sigma}^2_T = T^{-1} \int_0^T \sigma^2(X_t) \, dt$; it is an efficient estimator of $\|\sigma\|^2_{L^2(f)}$. An almost unbiased estimate of $|\varphi_{\sigma^2 f}(\lambda)|^2$ is then $|\hat{\varphi}_T(\lambda)|^2 - 4\hat{\sigma}^2_T/(T\lambda^2)$. The empirical characteristic function in this case has the form $\hat{\varphi}_T(\lambda) = T^{-1} \int_0^T e^{i\lambda X_t} \sigma^2(X_t) \, dt$. Accordingly, the functional $l_T$ is defined in this case as

$$l_T(h) = T \int_{\mathbb{R}} \lambda^2 (h^2(\lambda) - 2h(\lambda)) |\hat{\varphi}_T(\lambda)|^2 \, d\lambda + 8\hat{\sigma}^2_T \int_{\mathbb{R}} h(\lambda) \, d\lambda.$$

The remaining steps of the construction of the adaptive procedure do not need any modification. That is, we define $\tilde{h}_T$ by $l_T(\tilde{h}_T) = \min_{h \in \mathcal{H}^N_T} l_T(h)$. Next we apply to $\tilde{h}_T$ the inverse Fourier transform in order to define the kernel

$$\tilde{K}_T(x) = \frac{1}{2\pi} \int_{\mathbb{R}} \tilde{h}_T(\lambda) \cos(\lambda x) \, d\lambda.$$

Further, we choose another kernel function $Q(\cdot)$, which is positive, differentiable, symmetric, supported in $[-1, 1]$ and with $\int Q(u) \, du = 1$. Finally, we set $\varepsilon_T = e^{\sqrt{\log T}}$, $\ell_T = 1/\log T$ and define the estimator

$$\tilde{S}_T(x) = \frac{\int_0^T \sigma^2(X_t) \tilde{K}'_T(x - X_t) \, dt}{2\sqrt{T} \int_0^T Q(\sqrt{T}(x - X_t)) \, dt + 2\sqrt{T} \varepsilon_T e^{-\ell_T |x|}}.$$

The only thing that changes in Theorem 1 is the limiting constant. In this case the choice of a specific weight function does not allow one to obtain a limiting bound independent of $S_0$. The constant that we obtain is

$$P(S_0, \sigma, k, R) = P(k, R) \left( \int_{\mathbb{R}} \sigma^2(x) f_0(x) \, dx \right)^{2k/(2k+1)}.$$

4.3. *What happens if the diffusion is not ergodic?*. Note first that a diffusion, like any Markov process, can be positively recurrent, null recurrent or transient. Our method of adaptation, as well as the other methods suggested in the statistical literature for estimating adaptively the drift function, uses heavily the fact that the variance of the stochastic component in the risk decomposition is of order $1/(T\alpha_T)$, where $\alpha_T$ is the bandwidth or a smoothing parameter. This condition, as can be derived from the asymptotic equivalence result proven in [7], is not satisfied in the case of a null recurrent diffusion. The variance in that case is of order $1/(\sqrt{T}\alpha_T)$ and, consequently, the rates of convergence of drift estimators are significantly worse.

As to transient diffusions, even the simple feature of consistency fails for any estimator, since the amount of information concerning the value of $S$ at a point $x$ contained in the observed path $X^T$ does not increase when $T$ increases to infinity. That is the main reason for separating the ergodic case from the others.

4.4. *Sufficient conditions for* C3. A wide class of drift functions $S$ satisfying C3 is the set of all bounded functions: it is proven in [16] that there exist two positive constants $c_1$ and $c_2$ such that $p_t(S, x, y) \leq c_1 t^{-1/2} e^{-c_2|x-y|^2}$.

In the case when $X_t$ follows the invariant law, condition C3 is satisfied (with any $q < 2$ and any $\kappa > 0$) if $S$ is differentiable, satisfies condition C1 and $S(x)^2 + S'(x) > c$ for some constant $c \in \mathbb{R}$. Indeed, formula (7) on page 95 in [18] implies

$$p_t(S, x, y) \leq e^{-ct/2}\phi\left(\frac{x-y}{t}\right)\sqrt{f_S(y)/f_S(x)} \leq Cf_S(x)^{-1/2},$$

where $\phi(\cdot)$ stands for the density of the standard normal law. Therefore,

$$\int_{\mathbb{R}} \sup_y p_t^q(S, x, y)f_S(x)\,dx \leq C\int_{\mathbb{R}} f_S^{1-q/2}(x)\,dx < \infty,$$

since, under condition C1, the function $f_S$ decreases exponentially fast.

## 5. Proofs.

5.1. *An auxiliary result.* We start with a proposition reducing the study of the performance of drift estimators to that of the invariant density and its derivative estimators. From now on we suppose for simplicity that $\sigma \equiv 1$. Suppose now that $\bar{f}_T(\cdot)$ and $\bar{f}_T^{(1)}(\cdot)$ are estimators of the invariant density $f_S(\cdot)$ and its derivative $f_S'(\cdot)$, satisfying the conditions

E1. There exist $C_1, \gamma > 0$ such that $\mathbf{E}_S[\bar{f}_T(x) - f_S(x)]^{2p} \leq C_1^p T^{-p} e^{-\gamma|x|}$ for any $x \in \mathbb{R}$, $p \geq 1$ and $S \in \widehat{\Sigma}_\delta$.

E2. There exist two positive numbers $C_2$ and $q$ such that $\mathbf{E}_S[\bar{f}_T^{(1)}(x)^4] \leq C_2 T^q$, for any $x \in \mathbb{R}$ and $S \in \widetilde{\Sigma}_\delta$.

E3. The estimator $\bar{f}_T^{(1)}(\cdot)$ is asymptotically efficient, that is,

$$\lim_{\delta \to 0} \lim_{T \to \infty} \sup_{S \in \widetilde{\Sigma}_\delta} T^{2k/(2k+1)} \mathbf{E}_S \int_{\mathbb{R}} (\bar{f}_T^{(1)}(x) - f_S'(x))^2 \, dx = 4P(k, R).$$

Following some heuristics related to the identity $S(x) = f_S'(x)/2f_S(x)$ and presented in Section 2.1, we define the estimator of $S(x)$ as

$$(12) \qquad \hat{S}_T(x) = \frac{\bar{f}_T^{(1)}(x)}{2\bar{f}_T(x) + 2T^{-1/2}\varepsilon_T \, e^{-\ell_T |x|}},$$

where $\varepsilon_T = T^{1/\sqrt{\log T}} = e^{\sqrt{\log T}}$ and $\ell_T = (\log T)^{-1}$.

PROPOSITION 1.  *If conditions* E1–E3 *are fulfilled and* $S_0 \in \Sigma(k)$ *satisfies* C4, *then we have*

$$\lim_{\delta \to 0} \lim_{T \to \infty} \sup_{S \in \widetilde{\Sigma}_\delta} T^{2k/(2k+1)} R_T(\hat{S}_T, S) = P(k, R);$$

*that is, the estimator* $\hat{S}_T$ *is asymptotically minimax.*

PROOF.  The proof of this result relies on the Markov inequality and the exponential inequalities proven in Lemma 4 of [6]. It is quite similar to the proofs of Theorems 4 and 5 of [6] and therefore will be omitted here. For more details, we refer the reader to Theorem 6 of [5].  □

5.2. *Proof of Theorem* 1.  In the sequel the letters $C$ and $D$ stand for generic constants; the notation $\|h\|$ is used for the $L^2(\mathbb{R}, dx)$-norm of a function $h$. We assume that the initial value $\xi$ follows the invariant law; thus, the process $X$ is stationary in the strict sense.

Note that the estimator $\tilde{S}_T$ defined in Section 2.3 is of the form (12) with $\bar{f}_T^{(1)}(x) = T^{-1} \int_0^T \tilde{K}_T'(x - X_t) \, dt$ and $\bar{f}_T(x) = T^{-1/2} \int_0^T Q(\sqrt{T}(x - X_t)) \, dt$. Therefore, it suffices to verify that conditions E1–E3 hold. It is easy to see that Lemmas 4 and 5 and the arguments of Section 4.3 in [6] yield E1.

Condition E2 states that the fourth moment of $\bar{f}_T^{(1)}(x)$ is bounded in $x$ and increases in $T$ at most as a polynomial. This condition is evidently fulfilled, since

$$\tilde{K}_T'(x) = -(2\pi)^{-1} \int_{\mathbb{R}} \lambda \tilde{h}_T(\lambda) \sin(\lambda x) \, d\lambda$$

and the integrand above is bounded in absolute value by $|\lambda| \mathbb{1}_{[-T^{1/3}, T^{1/3}]}(\lambda)$ (recall that $\tilde{h}_T$ is supported by the interval $[-\tilde{\alpha}_T^{-1}, \tilde{\alpha}_T^{-1}]$, which is a subset of $[-T^{1/3}, T^{1/3}]$). Therefore, $\mathbf{E}_S[\bar{f}_T^{(1)}(x)^4]$ is bounded by $CT^{8/3}$.

It remains to verify E3, which is the most important part of the proof. For any estimator $\hat{f}_T^{(1)}(\cdot)$ of $f_S'(\cdot)$, we define the risk $r_T(\hat{f}_T^{(1)}, f_S')$ as the mean integrated squared error, that is, $r_T(\hat{f}_T^{(1)}, f_S') = \mathbf{E}_S \|\hat{f}_T^{(1)} - f_S'\|^2$. Due to the Plancherel identity, we have

$$r_T(\bar{f}_T^{(1)}, f_S') = \frac{1}{2\pi} \mathbf{E}_S \int_{\mathbb{R}} |\varphi_{\bar{f}_T^{(1)}}(\lambda) - \varphi_{f'}(\lambda)|^2 \, d\lambda$$

$$= \frac{1}{2\pi} \mathbf{E}_S \int_{\mathbb{R}} |\hat{\varphi}_T(\lambda) \varphi_{\tilde{K}'}(\lambda) - \varphi_{f'}(\lambda)|^2 \, d\lambda,$$

where we have used the notation $\hat{\varphi}_T(\lambda) = T^{-1} \int_{\mathbb{R}} e^{i\lambda X_t} \, dt$ and the fact that the Fourier transform of the convolution of two functions is the product of the Fourier transforms of the functions. Now, due to the formula of the Fourier transform of a derivative, we have

$$r_T(\bar{f}_T^{(1)}, f_S') = \frac{1}{2\pi} \mathbf{E}_S \int_{\mathbb{R}} |\lambda|^2 |[(\hat{\varphi}_T - \varphi_f)(\lambda)] \tilde{h}_T(\lambda) - \varphi_f(\lambda)(1 - \tilde{h}_T(\lambda))|^2 \, d\lambda.$$

This latter form of the risk is convenient since the term $\hat{\varphi}_T(\lambda) - \varphi_f(\lambda)$ is unbiased. Unfortunately, the randomness of the function $\tilde{h}_T$ does not allow us to apply to the risk $r_T$ the standard bias-variance decomposition. To bound this risk, more careful treatment of the main part of the stochastic component is required.

LEMMA 1. *For any $\lambda \in \mathbb{R}$, we have*

$$\lambda(\hat{\varphi}_T(\lambda) - \varphi_f(\lambda)) = 2i\zeta_T(\lambda) + T^{-1/2} m_S(\lambda, X^T),$$

*where $\zeta_T(\lambda) = T^{-1/2} \int_0^T e^{i\lambda X_t} \, dW_t$ and $m_S(\lambda, X^T)$ is a measurable function taking complex values such that, for sufficiently small values of $\delta > 0$, $\sup_{S \in \Sigma_\delta} \int_0^T \mathbf{E}_S |m_S(\lambda, X^T)|^2 \, d\lambda < C$.*

From now on, for two functions of $T$, say, $a_T$ and $b_T$, we write $a_T \sim b_T$ if the function $a_T/b_T$ tends to one as $T \to \infty$ uniformly in all the parameters entering in the definitions of these functions (in particular, uniformly in $S \in \Sigma_\delta$, for sufficiently small values of $\delta$). Using Lemma 1 and the fact that $\alpha_i \geq T^{-1/3}$, one can show that

$$R_T(\bar{f}_T^{(1)}, f_S') \sim \frac{1}{2\pi} \mathbf{E}_S \int_{\mathbb{R}} |2iT^{-1/2}\zeta_T(\lambda)\tilde{h}_T(\lambda) + \lambda\varphi_f(\lambda)(1 - \tilde{h}_T(\lambda))|^2 \, d\lambda.$$

The last expression is obviously of the same order as

$$\frac{4}{2\pi T} \mathbf{E}_S \int_{\mathbb{R}} |\zeta_T(\lambda)\tilde{h}_T(\lambda)|^2 \, d\lambda + \frac{1}{2\pi} \mathbf{E}_S \int_{\mathbb{R}} (1 - \tilde{h}_T(\lambda))^2 |\lambda\varphi_f(\lambda)|^2 \, d\lambda$$

$$- \frac{4}{2\pi\sqrt{T}} \operatorname{Im} \mathbf{E}_S \int_{\mathbb{R}} \lambda\tilde{h}_T(\lambda)(1 - \tilde{h}_T(\lambda))\zeta_T(\lambda)\varphi_f(-\lambda) \, d\lambda,$$

where $\text{Im}\, z$ is the imaginary part of the complex number $z$. This relation can be rewritten as

$$(13) \qquad R_T(\bar{f}_T^{(1)}, f_S') \sim \frac{1}{2\pi T}(\mathbf{E}_S[\Delta_T(1, \tilde{h}_T, \varphi_f)] + A_1 - \text{Im}\, A_2),$$

where $A_1 = 4\mathbf{E}_S \int_{\mathbb{R}} \tilde{h}_T^2(\lambda)(|\zeta_T(\lambda)|^2 - 1)\, d\lambda$ and

$$A_2 = 4\sqrt{T}\, \mathbf{E}_S \int_{\mathbb{R}} \lambda \tilde{h}_T(\lambda)(1 - \tilde{h}_T(\lambda))\zeta_T(\lambda)\varphi_f(-\lambda)\, d\lambda.$$

Note also that, from the definition of the functional $l_T$, one gets

$$l_T(h) = \int_{\mathbb{R}} 8h(\lambda)\, d\lambda + T \int_{\mathbb{R}} (h^2(\lambda) - 2h(\lambda))|\lambda \varphi_f(\lambda)|^2\, d\lambda$$

$$- 2T\, \text{Re} \int_{\mathbb{R}} (h^2(\lambda) - 2h(\lambda))\lambda \varphi_f(-\lambda)(\hat{\varphi}_T(\lambda) - \varphi_f(\lambda))\, d\lambda$$

$$+ T \int_{\mathbb{R}} \lambda^2(h^2(\lambda) - 2h(\lambda))|\hat{\varphi}_T(\lambda) - \varphi_f(\lambda)|^2\, d\lambda.$$

Using once more Lemma 1, we get

$$l_T(h) \sim \Delta_T(1, h, \varphi_f) + 4\sqrt{T}\, \text{Im} \int_{\mathbb{R}} (h^2(\lambda) - 2h(\lambda))\lambda \varphi_f(-\lambda)\, \zeta_T(\lambda)\, d\lambda$$

$$+ 4 \int_{\mathbb{R}} (h^2(\lambda) - 2h(\lambda))(|\zeta_T(\lambda)|^2 - 1)\, d\lambda - T\|\varphi_{f'}\|^2,$$

uniformly in $h \in \mathcal{H}_T$. This relation yields

$$(14) \qquad \mathbf{E}_S[l_T(\tilde{h}_T)] \sim \mathbf{E}_S[\Delta_T(1, \tilde{h}_T, \varphi_f)] - T\|\varphi_{f'}\|^2 + \text{Im}\, A_3 + A_4,$$

with the notation $A_3 = 4\sqrt{T}\, \mathbf{E}_S \int_{\mathbb{R}} (\tilde{h}_T^2(\lambda) - 2\tilde{h}_T(\lambda))\lambda \varphi_f(-\lambda)\zeta_T(\lambda)\, d\lambda$ and

$$A_4 = 4\mathbf{E}_S \int_{\mathbb{R}} (\tilde{h}_T^2(\lambda) - 2\tilde{h}_T(\lambda))(|\zeta_T(\lambda)|^2 - 1)\, d\lambda.$$

We wish to show that the terms $A_1$–$A_4$ are asymptotically smaller than $\mathbf{E}_S[\Delta_T(1, \tilde{h}, \varphi_f)]$ when $T$ tends to infinity. This can be done using the two following lemmas, the proofs of which are deferred to Section 5.3.

LEMMA 2. *Let $h(\cdot, w)$ be a bounded random function which takes only $N$ different values $h_1, \ldots, h_N$. Then $\mathbf{E}_S| \int_{\mathbb{R}} h(\lambda)\zeta_T(\lambda)\, d\lambda| \leq C\sqrt{N\mathbf{E}_S\|h\|^2}$, where the constant $C$ depends only on $k, R, S_0$.*

As a consequence of this lemma, we obtain $A_2 \vee A_3 \leq C\sqrt{N\mathbf{E}_S[\Delta_T(1, \tilde{h}, \varphi_f)]}$.

LEMMA 3. *For any random function $h(\cdot, \omega)$ taking only $N$ different values $h_1(\cdot), \ldots, h_N(\cdot)$ such that $\|h_i\|^2 \leq T$, the following inequality holds:*

$$\sup_{S \in \widetilde{\Sigma}_\delta} \mathbf{E}_S\left[ \int_{\mathbb{R}} h(\lambda)(|\zeta_T(\lambda)|^2 - 1)\, d\lambda \right] \leq C\sqrt{N\varepsilon_T}\sqrt{\mathbf{E}_S\|h\|^2},$$

where $\varepsilon_T = T^{1/\sqrt{\log T}}$ and the constant $C$ depends only on $k, R, S_0$.

As a consequence of this lemma, we obtain the inequality

$$\left| \mathbf{E}_S \left[ \int_{\mathbb{R}} \tilde{h}_T^n(\lambda)(|\zeta_T(\lambda)|^2 - 1) \, d\lambda \right] \right| \leq C\sqrt{N\varepsilon_T} \sqrt{\mathbf{E}_S \|\tilde{h}_T\|^2}$$
$$\leq C\varepsilon_T \sqrt{\Delta_T(1, \tilde{h}_T, \varphi_f)},$$

for any integer $n > 0$. This inequality implies that $A_1 \leq C\varepsilon_T \sqrt{\Delta_T(1, \tilde{h}_T, \varphi_f)}$ and $A_4 \leq C\varepsilon_T \sqrt{\Delta_T(1, \tilde{h}_T, \varphi_f)}$. Now (13) and (14) can be rewritten in the form

$$R_T(\bar{f}_T^{(1)}, f_S') \leq \frac{1}{2\pi T}(\sqrt{\mathbf{E}_S[\Delta_T(1, \tilde{h}_T, \varphi_f)]} + C\varepsilon_T)^2,$$

$$\mathbf{E}_S[l_T(h_T^*)] \leq (\sqrt{\mathbf{E}_S[\Delta_T(1, h_T^*, \varphi_f)]} + C\varepsilon_T)^2 - T\|\varphi_{f'}\|^2,$$

$$\mathbf{E}_S[l_T(\tilde{h}_T)] \geq (\sqrt{\mathbf{E}_S[\Delta_T(1, \tilde{h}_T, \varphi_f)]} - C\varepsilon_T)^2 - T\|\varphi_{f'}\|^2.$$

Here $h_T^*(\lambda) = (1 - |\alpha_T^* \lambda|^{k+\rho_T})_+$. Taking into account the fact that $\tilde{h}_T$ minimizes the functional $l_T(\cdot)$ over $\mathcal{H}_T^N$, we get $\mathbf{E}_S[l_T(\tilde{h}_T)] \leq \min_{h \in \mathcal{H}_T^N} \mathbf{E}_S[l_T(h)]$. On the other hand, by arguments very similar to those of Lemma 5 in [4], one checks easily that $\min_{h \in \mathcal{H}_T^N} \mathbf{E}_S[l_T(h)] \sim \min_{h \in \mathcal{H}_T} \mathbf{E}_S[l_T(h)] \leq \mathbf{E}_S[l_T(h_T^*)]$. Combining all these estimates, we arrive at the inequality

$$R_T(\bar{f}_T^{(1)}, f_S') \leq \frac{1}{2\pi T}(\sqrt{\Delta_T(1, h_T^*, \varphi_f)} + C\varepsilon_T)^2 (1 + o_T(1)).$$

Therefore, the expression $T^{k/(2k+1)}\sqrt{R_T(\bar{f}_T^{(1)}, f_S')}$ is asymptotically bounded by $\sqrt{T^{-1/(2k+1)}\Delta_T(1, h_T^*, \varphi_f)}$, plus a residual term $CT^{1/\sqrt{\log T} - 1/(4k+2)}$. It is well known in the theory of minimax estimation that the supremum of the quantity $T^{-1/(2k+1)}\Delta_T(1, h_T^*, \varphi_f)$ over the Sobolev ball $\Sigma(k+1, 4R) = \{f : \|f^{(k+1)} - f_0^{(k+1)}\|^2 \leq 4R\}$ tends to the constant $4P(k, R)$. This completes the proof of the theorem, due to the inclusion $\tilde{\Sigma}(k, R) \subset \Sigma(k+1, 4R + o_\delta(1))$ (see the proof of Theorem 5 in [6]).

### 5.3. *Proofs of technical lemmas.*

PROOF OF LEMMA 1. First of all, note that $\mathbf{E}_S[\hat{\varphi}_T(\lambda)] = \varphi_{f_S}(\lambda)$. Now, taking into account the occupation times formula ([33], page 224) and the martingale representation of the local time estimator ([26], page 137), we obtain $\hat{\varphi}_T(\lambda) - \mathbf{E}_S[\hat{\varphi}_T(\lambda)] = T^{-1}(H_S(\lambda, X_T) - H_S(\lambda, X_0)) - T^{-1}\int_0^T g_S(\lambda, X_t) \, dW_t,$ where the functions $H_S$ and $g_S$ are defined as

$$g_S(\lambda, u) = 2 \int_{\mathbb{R}} e^{i\lambda x} f_S(x) \left( \frac{\mathbb{1}_{\{u > x\}} - F_S(u)}{f_S(u)} \right) dx,$$

$$H_S(\lambda, u) = 2 \int_{\mathbb{R}} e^{i\lambda x} f_S(x) \left[ \int_0^u \left( \frac{\mathbb{1}_{\{v > x\}} - F_S(v)}{f_S(v)} \right) dv \right] dx.$$

The integration by parts formula yields

$$i\lambda g_S(\lambda, u) = 2 e^{i\lambda u} - \tilde{g}_S(\lambda, u),$$

where

$$\tilde{g}_S(\lambda, u) = 2 \int_{\mathbb{R}} e^{i\lambda x} f_S'(x) \left( \frac{\mathbb{1}_{\{u > x\}} - F_S(u)}{f_S(u)} \right) dx.$$

It implies that $\lambda(\hat{\varphi}_T(\lambda) - \varphi_f(\lambda)) = T^{-1} 2i \int_0^T e^{i\lambda X_t} \, dW_t + T^{-1/2} m_S(\lambda, X^T)$ with $\sqrt{T} m_S(\lambda, X^T) = \lambda(H_S(\lambda, X_T) - H_S(\lambda, X_0)) + i \int_0^T \tilde{g}_S(\lambda, X_t) \, dW_t$. In the same way one can prove that $|i\lambda H_S(\lambda, u)| = |\int_0^u i\lambda g_S(\lambda, v) \, dv| \leq 2|u| + |\int_0^u \tilde{g}_S(\lambda, v) \, dv|$. Using the Plancherel identity and Lemma 4 from [6], we get

$$\int_{\mathbb{R}} \mathbf{E}_S |\tilde{g}_S(\lambda, \xi)|^2 \, d\lambda = 2\pi \int_{\mathbb{R}} f_S'(x)^2 \mathbf{E}_S \left[ \frac{\mathbb{1}_{\{\xi > x\}} - F_S(\xi)}{f_S(\xi)} \right]^2 dx \leq C,$$

where $C$ is a constant independent of $S \in \Sigma_\delta$. Similarly, one checks that $\int_{\mathbb{R}} \mathbf{E}_S |\int_0^\xi \tilde{g}_S(\lambda, v) \, dv|^2 \, d\lambda \leq C$. Thus, we have

$$\int_0^T \mathbf{E}_S |m_S(\lambda, X^T)|^2 \, d\lambda \leq \int_0^T \left( \frac{8}{T} \mathbf{E}_S |i\lambda H_S(\lambda, \xi)|^2 + 2\mathbf{E}_S |g_S(\lambda, \xi)|^2 \right) d\lambda$$

$$\leq 16 \mathbf{E}_S[\xi^2] + 16 \int_{\mathbb{R}} \mathbf{E}_S \left| \int_0^\xi \tilde{g}_S(\lambda, v) \, dv \right|^2 d\lambda$$

$$+ 2 \int_{\mathbb{R}} \mathbf{E}_S |\tilde{g}_S(\lambda, \xi)|^2 \, d\lambda < C,$$

and the assertion of the lemma follows.  $\square$

Proof of Lemma 2.   Let us denote $\xi_h = \|h\|^{-1} \int_{\mathbb{R}} h \zeta_T$. It is evident that

$$\mathbf{E}_S \left| \int_{\mathbb{R}} h(\lambda) \zeta_T(\lambda) \, d\lambda \right| = \mathbf{E}_S(\|h\| \cdot |\xi_h|) \leq \sqrt{\mathbf{E}_S \|h\|^2 \mathbf{E}_S |\xi_h|^2}$$

$$\leq \left[ \mathbf{E}_S \|h\|^2 \sum_{i=1}^N \mathbf{E}_S |\xi_{h_i}|^2 \right]^{1/2}.$$

Now, taking into account the explicit form of $\zeta_T$, we have

$$\mathbf{E}_S |\xi_{h_i}|^2 = \frac{1}{T \|h_i\|^2} \mathbf{E}_S \left| \int_0^T \int_{\mathbb{R}} e^{i\lambda X_t} h_i(\lambda) \, d\lambda \, dW_t \right|^2$$

$$= \frac{1}{\|h_i\|^2} \mathbf{E}_S \left| \int_{\mathbb{R}} e^{i\lambda\xi} h_i(\lambda)\, d\lambda \right|^2 = \frac{1}{\|h_i\|^2} \mathbf{E}_S |\varphi_{h_i}(\xi)|^2$$

$$\leq \frac{C}{\|h_i\|^2} \left[ \int_{\mathbb{R}} |\varphi_{h_i}(x)|^2\, dx \right] = C,$$

where $C = \sup_{S \in \Sigma_\delta} \sup_{x \in \mathbb{R}} f_S(x)$. This completes the proof of Lemma 2. $\square$

PROOF OF LEMMA 3. The Itô formula implies that, for any continuous martingale $\mathcal{M}_t$, we have $\mathcal{M}_T^2 - \mathcal{M}_0^2 = 2\int_0^T \mathcal{M}_t\, d\mathcal{M}_t + \langle\mathcal{M}\rangle_T$, where $\langle\mathcal{M}\rangle_t$ is the quadratic variation of the martingale $\mathcal{M}_t$. Applying this formula, we get $T|\zeta_T(\lambda)|^2 = 2\int_0^T Y_t(\lambda)\, dW_t + T$, where we have used the notation $Y_t(\lambda) = \operatorname{Re} e^{i\lambda X_t} \int_0^t e^{-i\lambda X_u}\, dW_u = \operatorname{Re} e^{i\lambda X_t} \sqrt{t}\, \zeta_t(\lambda)$. Changing the order of the integrals and using the Itô isometry, we get

$$U_S(h) := \mathbf{E}_S \left| \int_{\mathbb{R}} h(\lambda)(|\zeta_T(\lambda)|^2 - 1)\, d\lambda \right|^2$$

$$= \frac{4}{T^2} \mathbf{E}_S \left[ \int_0^T \left| \int_{\mathbb{R}} Y_t(\lambda) h(\lambda)\, d\lambda \right|\, dW_t \right]^2$$

$$= \frac{4}{T^2} \int_0^T \mathbf{E}_S \left| \int_{\mathbb{R}} h(\lambda) Y_t(\lambda)\, d\lambda \right|^2\, dt$$

$$\leq \frac{4}{T^2} \int_0^T t\mathbf{E}_S \left| \int_{\mathbb{R}} e^{i\lambda X_t} \zeta_t(\lambda)\, h(\lambda)\, d\lambda \right|^2\, dt.$$

We apply now the same method as in the proof of the first lemma. Let us introduce $\xi_h = \|h\|^{-1} \int_{\mathbb{R}} h(\lambda)(|\zeta_T(\lambda)|^2 - 1)\, d\lambda$. The Cauchy–Schwarz inequality yields

$$\mathbf{E}_S \left[ \int_{\mathbb{R}} h(\lambda)(|\zeta_T(\lambda)|^2 - 1)\, d\lambda \right] \leq \sum_{i=1}^N \sqrt{\mathbf{E}_S \|h\|^2} \sqrt{\mathbf{E}_S |\xi_{h_i}|^2}.$$

It remains to carry out the inequality $\mathbf{E}_S |\xi_{h_i}|^2 = \|h_i\|^{-2} U_S(h_i) \leq C$, where $C$ is a constant depending only on $k, R, S_0$. The proof of this inequality will be divided into three steps.

*First step.* Suppose that $T > \kappa$, where $\kappa$ is the positive number defined by condition C3. Then we have the inequalities

$$\int_0^\kappa t\mathbf{E}_S \left| \int_{\mathbb{R}} e^{i\lambda X_t} \zeta_t(\lambda)\, h_i(\lambda)\, d\lambda \right|^2\, dt \leq \kappa^2 \left( \int_{\mathbb{R}} h_i(\lambda)\, d\lambda \right)^2 \leq 4\kappa^2 \alpha_i^{-2}$$

$$\leq 12\kappa^2 T \|h_i\|^2,$$

$$\mathbf{E}_S \left| \int_{\mathbb{R}} e^{i\lambda X_t} h_i(\lambda) \int_{t-\kappa}^t e^{i\lambda X_u}\, dW_u\, d\lambda \right|^2 \leq \mathbf{E}_S \left[ \int_{-\alpha_i^{-1}}^{\alpha_i^{-1}} \left| \int_{t-\kappa}^t e^{i\lambda X_u}\, dW_u \right|\, d\lambda \right]^2$$

$$\leq 2\alpha_i^{-1} \int_{-\alpha_i^{-1}}^{\alpha_i^{-1}} \mathbf{E}_S \Big| \int_{t-\kappa}^{t} e^{i\lambda X_u} \, dW_u \Big|^2 \, d\lambda$$

$$\leq 4\alpha_i^{-2}\kappa \leq 12\kappa T \|h_i\|^2.$$

These inequalities imply that

$$\frac{U_S(h_i)}{\|h_i\|^2} \leq C + \frac{4}{\|h_i\|^2 T^2} \int_{\kappa}^{T} \mathbf{E}_S \Big| \int_{\mathbb{R}} e^{i\lambda X_t} h_i(\lambda) \int_{0}^{t-\kappa} e^{i\lambda X_u} \, dW_u \, d\lambda \Big|^2 \, dt$$

$$= C + \frac{4\tilde{U}_S(h_i)}{\|h_i\|^2 T^2}.$$

We want to prove now that $\tilde{U}_S(h_i) \leq C\varepsilon_T T^2 \|h_i\|^2$ for a constant $C$. Recall that $X^u$ denotes the trajectory of $X$ between 0 and $u$. Since the random variable $\eta(X^{t-\kappa}, \lambda) = \int_{0}^{t-\kappa} e^{i\lambda X_u} \, dW_u$ is $\mathscr{F}_{t-\kappa}$-measurable and the law $\mathscr{L}(X_t|\mathscr{F}_{t-\kappa}) = \mathscr{L}(X_t|X_{t-\kappa})$, we have

$$\tilde{U}_S(h_i) = \int_{\kappa}^{T} \mathbf{E}_S \Big[ \int_{\mathbb{R}} \Big| \int_{\mathbb{R}} e^{i\lambda y} \eta(X^{t-\kappa}, \lambda) h_i(\lambda) \, d\lambda \Big|^2 p_\kappa(S, X_{t-\kappa}, y) \, dy \Big] \, dt.$$

Let us denote $Q(X^{t-\kappa}, y) = |\int_{\mathbb{R}} e^{i\lambda y} \eta(X^{t-\kappa}, \lambda) \, h_i(\lambda) \, d\lambda|$.

*Second step.* To simplify the notation we suppose that $\kappa = 1$. Let us prove now that, for any function $Q(\cdot)$, the following inequality holds:

$$\begin{aligned}
(15) \qquad \int_{\mathbb{R}} Q^2(y) p_1(S, x, y) \, dy \leq {} & \varepsilon_T \sup_{y \in \mathbb{R}} p_1(S_0, x, y) \int_{\mathbb{R}} Q^2(y) \, dy \\
& + \frac{1}{T^2} \sqrt{\int_{\mathbb{R}} Q^4(y) p_1(S, x, y) \, dy},
\end{aligned}$$

for any $S \in V_\delta(S_0)$ with $\delta \leq 0.2$. Indeed, if we denote by $\mathbf{E}_S^x$ the mathematical expectation with respect to the measure induced by the solution of (1) with deterministic initial value $X_0 = x$, then $\int_{\mathbb{R}} Q^2(y) \, p_1(S, x, y) \, dy = \mathbf{E}_S^x[Q^2(X_1)] = \mathbf{E}_{S_0}^x[Q^2(X_1)L(S, X^1)]$, where

$$L_1(S, X^1) = \exp\Big\{ \int_{0}^{1} (S(X_u) - S_0(X_u)) \, dW_u - \tfrac{1}{2} \int_{0}^{1} (S(X_u) - S_0(X_u))^2 \, du \Big\}$$

is the likelihood ratio and $X^1 = (X_t, 0 \leq t \leq 1)$ denotes the trajectory of the process $X$ up to the time 1. Further, using the Cauchy–Schwarz inequality, we get

$$\mathbf{E}_S^x[Q^2(X_1)] \leq \varepsilon_T \mathbf{E}_{S_0}^x[Q^2(X_1)] + \sqrt{\mathbf{E}_S^x[Q^4(X_1)]} \sqrt{\mathbf{P}_S^x(L(S, X^1) > \varepsilon_T)}.$$

Note now that, for any $n > 0$, we have

$$L(S, X^1)^n = \exp\Big\{ M_1 - \frac{1}{2}\langle M \rangle_1 + \frac{n^2 - n}{2} \int_{0}^{1} (S(X_u) - S_0(X_u))^2 \, du \Big\}$$

$$\leq e^{M_1 - (1/2)\langle M \rangle_1} \, e^{n^2 \delta^2},$$

where $M_t = n \int_0^t (S - S_0)(X_u)\, dW_u$ is a local martingale. Therefore, taking $n = 5\sqrt{\log T}$ and applying the Chebyshev inequality, we get $\mathbf{P}_S^x(L(S, X^1) > \varepsilon_T) \le \varepsilon_T^{-5\sqrt{\log T}} e^{25\delta^2 \log T} \le T^{-4}$, which leads to (15).

*Third step.* Inequality (15) yields $\tilde{U}_S(h_i) \le \int_1^T (D_1(t) + T^{-2}\sqrt{D_2(t)})\, dt$, where we have used the abbreviations

$$D_1(t) = \varepsilon_T \mathbf{E}_S \left[ \sup_{y \in \mathbb{R}} p_1(S_0, X_{t-1}, y) \int_{\mathbb{R}} \left| \int_{\mathbb{R}} e^{i\lambda y} h_i(\lambda) \eta(X^{t-1}, \lambda)\, d\lambda \right|^2 dy \right]$$

and $D_2(t) = \mathbf{E}_S \int_{\mathbb{R}} Q^4(X^{t-1}, y) p_1(S, X_{t-1}, y)\, dy = \mathbf{E}_S[Q^4(X^{t-1}, X_t)]$. The term $D_1$ can be evaluated via the Plancherel identity and the Hölder inequality:

$$D_1(t) = 2\varepsilon_T \pi \mathbf{E}_S \left[ \sup_y p_1(S_0, X_{t-1}, y) \int_{\mathbb{R}} |h_i(\lambda)\eta(X^{t-1}, \lambda)|^2\, d\lambda \right]$$

$$\le C_q \varepsilon_T (t-1) \left( \mathbf{E}_S \left[ \sup_y p_1(S_0, X_0, y)^q \right] \right)^{1/q} \|h_i\|^2,$$

in view of the estimate $\mathbf{E}[|\eta(X^{t-1}, \lambda)|^{2s}] \le C_s(t-1)^s$, for any $s > 0$, which follows from the BDG inequality ([33], Theorem IV.4.1). Using condition C3 and the [uniform in $S \in V_\delta(S_0)$] boundedness of $f_S/f_{S_0}$, we get $D_1(t) \le C\varepsilon_T(t-1)\|h_i\|^2$.

One easily checks that $|Q(X^{t-1}, X_t)|^2 \le \int_{|\alpha_i \lambda| < 1} |\eta(X^{t-1}, \lambda)|^2\, d\lambda \|h_i\|^2$ and, hence, $D_2(t) \le Ct^2\|h_i\|^8 \le Ct^2 T^2\|h_i\|^4$. Combining these estimates, one obtains $\tilde{U}_S(h_i) \le CT^2\|h_i\|^2$, which completes the proof. $\quad\square$

LABORATOIRE DE PROBABILITÉS
UNIVERSITÉ PARIS 6
BOÎTE COURRIER 188
75252 PARIS CEDEX 05
FRANCE
E-MAIL: dalalyan@ccr.jussieu.fr